\documentclass[10pt]{amsart}
\usepackage{amsthm,amsmath,amssymb,amsfonts}
\usepackage[american]{babel}
\usepackage{cite}
\usepackage{url}

\newtheorem{theorem}{Theorem}
\newtheorem*{theorem*}{Theorem}

\theoremstyle{definition}

\newtheorem{example}[theorem]{Example}
\newtheorem*{example*}{Example}

\theoremstyle{remark}
\newtheorem{remark}[theorem]{Remark}

\newtheorem*{remark*}{Remark}

\begin{document}


\newcommand{\abs}[1]{\lvert#1\rvert}
\newcommand{\nt}{\trianglelefteq}
\newcommand{\NU}[1]{\mbox{\rm V}(#1)}
\newcommand{\cen}[2]{\mbox{\rm C}_{#1}(#2)}
\newcommand{\zen}[1]{\mbox{\rm Z}(#1)}
\newcommand{\lara}[1]{\langle{#1}\rangle}
\newcommand{\ZZ}{\mbox{$\mathbb{Z}$}}
\newcommand{\QQ}{\mbox{$\mathbb{Q}$}}
\newcommand{\sQQ}{\mbox{$\scriptstyle\mathbb{Q}$}}
\newcommand{\paug}[2]{\varepsilon_{#1}(#2)}
\newcommand{\eaug}[1]{\varepsilon_{{\tt #1}}}

\newenvironment{mydescription}[1]
{\begin{list}{}
{\renewcommand{\makelabel}[1]{#1}
\settowidth{\labelwidth}{#1}
\settowidth{\leftmargin}{0pt}
\addtolength{\leftmargin}{\labelwidth}
\addtolength{\leftmargin}{-2\labelsep}
\addtolength{\leftmargin}{1pt}
}}{\end{list}}


\title[Zassenhaus conjecture for central extensions of $S_{5}$]
{Zassenhaus conjecture for \\
central extensions of $\boldsymbol{S}_{\boldsymbol{5}}$}

\author{Victor Bovdi}
\address{Institute of Mathematics, University of Debrecen,
4010 Debrecen, P.O.~Box 12, \mbox{Hungary};
Institute of Mathematics and Informatics,
College of Ny\'{i}regyh\'{a}za, S\'{o}st\'{o}i \'{u}t 31/B,
4410 Ny\'{i}regyh\'{a}za, Hungary}
\email{vbovdi@math.klte.hu}

\author{Martin Hertweck}
\address{Universit\"at Stuttgart, Fachbereich Mathematik, IGT,
Pfaffenwald\-ring 57, 70550 Stuttgart, Germany}
\email{hertweck@mathematik.uni-stuttgart.de}

\footnotetext{The research  was supported by OTKA T 037202, T
038059}



\subjclass[2000]{Primary 16S34, 16U60; Secondary 20C05}


\keywords{integral group ring, torsion unit, Zassenhaus conjecture}

\begin{abstract}
We confirm a conjecture of Zassenhaus about rational conjugacy of
torsion units in integral group rings for a covering group of the
symmetric group $S_{5}$ and for the general linear group $\text{GL}(2,5)$.
\end{abstract}

\maketitle

\section{Introduction}\label{intro}

The conjecture of the title states:
\begin{mydescription}{(ZC1)}
\item[(ZC1)]
For a finite group $G$, every torsion unit in $\ZZ G$ is conjugate to
an element of $\pm G$ in the units of $\QQ G$.
\end{mydescription}

It remains not only unsolved but also lacking in
plausible means of either proof or counter-example, at least for
non-solvable groups $G$.
The purpose of this note is to add two further groups to the small
list of non-solvable groups $G$ for which conjecture (ZC1) has been
verified (see \cite{LuPa:89,LuTr:91,DoJuPM:97,He:05b}).

\begin{example}\label{ex1}
The conjecture (ZC1) holds true for
the covering group $\widetilde{S}_{5}$ of the symmetric group
$S_{5}$ which contains a unique conjugacy class of involutions.
\end{example}

\begin{example}\label{ex2}
The conjecture (ZC1) holds true for
the general linear group $\text{GL}(2,5)$.
\end{example}

We remark that $\text{\rm PGL}(2,5)\cong S_{5}$
(see \cite[Kapitel~II, 6.14~Satz]{Hupp:67}).

The covering group $\widetilde{S}_{5}$ occurs as  Frobenius complement
in Frobenius groups (for the classification of Frobenius complements
see \cite{Pass:68}).
From already existing work \cite{DoJu:96,DoJuPM:97,JuPM:00} it follows
that Example~\ref{ex1} supplies the missing bit for the proof of:
\begin{theorem}
Let $G$ be a finite Frobenius group. Then each torsion unit
in $\ZZ G$ which is of prime power order is conjugate to
an element of $\pm G$ in the units of $\QQ G$.
\end{theorem}

The proofs are obtained by applying the Luthar--Passi method
\cite{LuPa:89,He:05b}. We dispose of the validity of (ZC1) for
$S_{5}$, established in \cite{LuTr:91} (see also
\cite[Section~4]{He:05b} for a proof using the Luthar--Passi
method). Below, we briefly recall this method, which uses the
character table and/or modular character tables in an automated
process suited for being done on a computer, the result being that
rational conjugacy of torsion units of a given order to group
elements is either proven or not, and if not, one gets at least
some information about partial augmentations (cf.\
\cite{BoHoKi:04,BoKo:06,BoJeKo:06}). We tried to prevent the
proofs from useless ballast and to present them in human-readable
format, rather than producing systems of inequalities and their
solutions which should be reasonably done on a computer (cf.
\cite{LAGUNA}).

Let $G$ be a finite group.
Recall that for a group ring element $u=\sum_{g\in G}a_{g}g$
(all $a_{g}\in\ZZ$), its partial augmentation
with respect to the conjugacy class $x^{G}$ of an element $x$ of $G$,
in the following denoted by $\paug{x}{u}$ or $\paug{x^{G}}{u}$,
is $\sum_{g\in x^{G}}a_{g}$.
When dealing with conjecture (ZC1), it suffices to consider units
of augmentation one which form a group denoted by $\NU{\ZZ G}$.

A few preliminary remarks on a torsion unit $u$ in $\NU{\ZZ G}$
seem to be appropriate.

The familiar result of Berman--Higman
(from \cite{Ber:55} and \cite[p.~27]{Hig:40})
asserts that if $\paug{z}{g}\neq 0$ for some $z$ in the center
of $G$, then $u=z$.

A practical criterion for $u$ being conjugate to an element of $G$ in
the units of $\QQ G$ is that all but one of the partial augmentations
of every power of $u$ must vanish (see \cite[Theorem~2.5]{MaRiSeWe:87}).
The converse is obvious.

The next remarks will be used repeatedly.

\begin{remark}\label{rem1}
Let $u$ be a torsion unit in $\NU{\ZZ G}$, let $N\nt G$ and set $\bar{G}=G/N$.
We write $\bar{u}$ for the image of
$u$ under the natural map $\ZZ G\rightarrow\ZZ\bar{G}$.
Since any conjugacy class of $G$ maps
onto a conjugacy class of $\bar{G}$, we have that for any $x\in G$,
the partial augmentation $\paug{\bar{x}}{\bar{u}}$ is the sum of
the partial augmentations $\paug{g^{G}}{u}$ with $g\in G$ such that
$\bar{g}$ is conjugate to $\bar{x}$ in $\bar{G}$.

Now suppose that $N$ is a central subgroup of $G$, and that $\bar{u}=1$.
Then $u\in N$. Indeed, $1=\paug{1}{\bar{u}} = \sum_{n\in N}\paug{n}{u}$,
so $u$ has a central group element in its support and the
Berman--Higman result applies.
\end{remark}

\begin{remark}\label{rem2}
Let $u$ be a torsion unit in $\NU{\ZZ G}$. Then $g\in G$ and
$\paug{g}{u}\neq 0$ implies that the order of $g$ divides the order of $u$.
Indeed, it is well known that then prime divisors of the order
of $g$ divide the order of $u$ (see \cite[Theorem~2.7]{MaRiSeWe:87},
as well as \cite[Lemma~2.8]{He:05a} for an alternative
proof). Further, it was observed in \cite[Lemma~5.6]{He:05b} that
the orders of the $p$-parts of $g$ cannot exceed those of $u$.
\end{remark}

We briefly recall the Luthar--Passi method. Let $u\in\NU{\ZZ G}$.
Suppose that $u^{n}=1$ for some natural number $n$ and let $\zeta$
be a primitive complex $n$-th root of unity. Let $\chi$ be the character
afforded by a complex representation $D$ of $G$, and write
$\mu(\xi,u,\chi)$ for the multiplicity of an $n$-th root of unity $\xi$
as an eigenvalue of the matrix $D(u)$. Then
(cf.\ \cite{LuPa:89}, \cite[Section~3]{He:05b}):
\[ \mu(\xi,u,\chi) = \frac{1}{n}\sum_{d\mid n}
\text{\rm Tr}_{\sQQ(\zeta^{d})/\sQQ}(\chi(u^{d})\xi^{-d}). \]

When trying to show that $u$ is rationally conjugate to an element of $G$,
one may assume---by induction on the order of $u$---that the values of the
summands for $d\neq 1$ are ``known.'' The summand for $d=1$ can be
written as
\[ \frac{1}{n}\sum_{g^{G}}\paug{g}{u}
\text{\rm Tr}_{\sQQ(\zeta)/\sQQ}(\chi(g)\xi^{-1}), \]
a linear combination
of the $\paug{g}{u}$ with ``known'' coefficients. Note that the
$\mu(\xi,u,\chi)$ are non-negative integers, bounded above by $\chi(1)$.
Thus in some sense, there are linear inequalities in the partial augmentations
of $u$ which impose constraints on them. Trying to make use of these
inequalities is now understood as being the Luthar--Passi method.

A modular version of this method (see \cite[Section~3]{He:05b} for
details) can be derived from the following observation in the very same way
as the original (complex) version is derived from the (obvious) fact that
$\chi(u)=\sum_{g^{G}}\paug{g}{u}\chi(g)$.
Suppose that $p$ is a rational prime which does not divide the order of $u$
(i.e., $u$ is a $p$-regular torsion unit).
Then for every Brauer character $\varphi$ of $G$ (relative to $p$) we have
(see \cite[Theorem~2.2]{He:05b}):
\[ \varphi(u)=\sum_{\substack{g^{G}:\; g\text{ is}\\ \text{$p$-regular}}}
\paug{g}{u}\varphi(g). \]
Thereby, the domain of $\varphi$ is naturally extended to the set of
$p$-regular torsion units in $\ZZ G$.

\section{A covering group of $S_{5}$}\label{Sec:S5cover}

A presentation of a covering group of $S_{n}$ is given by
\begin{align*}
\widetilde{S}_{n}=\langle\,
g_{1}, \dotsc, g_{n-1}, z
\; | \;&
g_{i}^{2}=(g_{j}g_{j+1})^{3}=(g_{k}g_{l})^{2}=z,\,
z^{2}=[z,g_{i}]=1 \\
& \text{for } 1\leq i\leq n-1, \, 1\leq j\leq n-2, \, k\leq l-2\leq n-3
\, \rangle.
\end{align*}
Recent results in the representation theory of the covering groups of
symmetric groups can be found in Bessenrodt's survey article \cite{Bess:94}.
We merely remark that the complex spin characters of $\widetilde{S}_{n}$,
i.e., those characters which are not characters of $S_{n}$,
were determined by Schur \cite{Schu:11}.

The group $\widetilde{S}_{5}$ has catalogue number $89$
in the Small Group Library in GAP \cite{GAP4}
(the other covering group of $S_{5}$ has number $90$).
The spin characters of $\widetilde{S}_{5}$ as produced by GAP
are shown in Table~\ref{T1} (dots indicate zeros).

\begin{table}[h]
\[
\begin{array}{c}
\begin{array}{r|rrrrrrrrrrrr} \hline
& {\tt 1a} & {\tt 5a} & {\tt 4a} & {\tt 2a} & {\tt 10a} & {\tt 6a}
& {\tt 3a} & {\tt 8a} & {\tt 8b} & {\tt 4b} & {\tt 12a} & {\tt 12b}
\rule[-7pt]{0pt}{20pt} \\ \hline
\chi_{5} & 4& -1&  . & -4&   1&  2& -2&  .&  .&  .&   .&  .
\rule[0pt]{0pt}{13pt} \\
\chi_{6} & 4& -1&  .& -4&   1& -1&  1&  .&  .&  .&   \beta&  -\beta\\
\chi_{7} & 4& -1&  .& -4&   1& -1&  1&  .&  .&  .&  -\beta&   \beta\\
\chi_{11}& 6&  1&  .& -6&  -1&  .&  .&  \alpha& -\alpha&  .&   .&   .\\
\chi_{12}& 6&  1&  .& -6&  -1&  .&  .& -\alpha&  \alpha&  .&   .&   .
\rule[-7pt]{0pt}{5pt} \\ \hline
\end{array} \\
\begin{array}{rl}
\text{Irrational entries:} &
\alpha=-\zeta_{8}^{}+\zeta_{8}^{3}=-\sqrt{2} \text{ where }
\zeta_{8}=\exp(2\pi i/8), \\
& \beta=\zeta_{12}^{7}-\zeta_{12}^{11}=-\sqrt{3} \text{ where }
\zeta_{12}=\exp(2\pi i/12).
\rule[6pt]{0pt}{5pt} \\
\end{array}
\rule[0pt]{0pt}{25pt}
\end{array} \]
\rule[0pt]{0pt}{5pt}
\caption{Spin characters of $\widetilde{S}_{5}$}\label{T1}
\end{table}

We turn to the proof that conjecture (ZC1)
holds true for $\widetilde{S}_{5}$.
Let $z$ be the central involution in $\widetilde{S}_{5}$.
Then we have a natural homomorphism
$\pi:\ZZ\widetilde{S}_{5}\rightarrow\ZZ\widetilde{S}_{5}/\lara{z}=\ZZ S_{5}$.
Let $u$ be a nontrivial torsion unit in $\NU{\ZZ\widetilde{S}_{5}}$.
We shall show that all but one of its partial augmentations vanish.
Since (ZC1) is true for $S_{5}$, the order of $\pi(u)$ agrees with the
order of an element of $S_{5}$, and it follows that the order of $u$
agrees with the order of an element of $\widetilde{S}_{5}$
(see Remark~\ref{rem1}). By the Berman--Higman result we can assume that
$\paug{1}{u}=0$ and  $\paug{z}{u}=0$.
Further, we can assume that the order of $u$ is even since otherwise
rational conjugacy of $u$ to an element of $G$ follows from
the validity of (ZC1) for $S_{5}$ and \cite[Theorem~2.2]{DoJu:96}.
Denote the partial augmentations of $u$ by
$\eaug{1a},\eaug{5a},\dotsc,\eaug{12b}$ (so that $\eaug{5a}$, for example,
denotes the partial augmentations of $u$ with respect to the conjugacy class
of elements of order $5$). So $\eaug{1a}=\eaug{2a}=0$.
Since all but one of the partial augmentations of $\pi(u)$, the image
of $u$ in $\ZZ S_{5}$, vanish, and a partial augmentation of $\pi(u)$ is
the sum of the partial augmentations of $u$ taken for classes which
fuse in $S_{5}$, we have
\begin{equation}\label{eq1}
\begin{split}
& \eaug{4a},\eaug{4b},\eaug{8a}+\eaug{8b},\eaug{3a}+\eaug{6a},
\eaug{5a}+\eaug{10a},\eaug{12a}+\eaug{12b} \in\{0,1\}, \\
& \abs{\eaug{4a}}+\abs{\eaug{4b}}+\abs{\eaug{8a}+\eaug{8b}}
+\abs{\eaug{3a}+\eaug{6a}}+\abs{\eaug{5a}+\eaug{10a}}+
\abs{\eaug{12a}+\eaug{12b}}=1.
\end{split}
\end{equation}

\subsection*{When $\boldsymbol{u}$ has order $\boldsymbol{2}$ or
$\boldsymbol{4}$.}
Then the partial augmentations of $u$ which are possibly nonzero are
$\eaug{4a}$, $\eaug{4b}$, $\eaug{8a}$ and $\eaug{8b}$
(Remark~\ref{rem2}). Thus
$\chi_{11}(u)=\alpha(\eaug{8a}-\eaug{8b})=-\sqrt{2}(\eaug{8a}-\eaug{8b})$.
Also, $\chi_{11}(u)$ is the sum of fourth roots of unity. Since
$\sqrt{2}\not\in\QQ(i)$ it follows that $\eaug{8a}-\eaug{8b}=0$.
Using $\eaug{8a}+\eaug{8b}\in\{0,1\}$ from (\ref{eq1}) we obtain
$\eaug{8a}=\eaug{8b}=0$. Now $\abs{\eaug{4a}}+\abs{\eaug{4b}}=1$ by
(\ref{eq1}), so all but one of the partial augmentations of $u$ vanish,
with either $\eaug{4a}=1$ or $\eaug{4b}=1$. It follows that $u$ is
rationally conjugate to an element of $G$ (necessarily of order $4$).

\subsection*{When $\boldsymbol{u}$ has order $\boldsymbol{6}$ or
$\boldsymbol{10}$.} Then $u^{3}=z$ or $u^{5}=z$, respectively
(Remark~\ref{rem1}), i.e.,
$zu$ is of order $3$ or $5$. Thus $zu$ is, as already noted, rationally
conjugate to an element of $G$, and hence the same holds for $u$ itself.

\subsection*{When $\boldsymbol{u}$ has order $\boldsymbol{12}$.}
Then the partial augmentations of $u$ which are possibly nonzero are
$\eaug{4a}$, $\eaug{4b}$, $\eaug{8a}$, $\eaug{8b}$,
$\eaug{3a}$, $\eaug{6a}$, $\eaug{12a}$ and $\eaug{12b}$.
The unit $\pi(u)$ has order $6$ (Remark~\ref{rem2}),
so $\eaug{12a}+\eaug{12b}=1$ and
$\eaug{4a}=\eaug{4b}=\eaug{8a}+\eaug{8b}=\eaug{3a}+\eaug{6a}=0$
by (\ref{eq1}). Now $\chi_{11}(u)=-\sqrt{2}(\eaug{8a}-\eaug{8b})$
but $\sqrt{2}\not\in\QQ(\zeta_{12})=\QQ(i,\zeta_{3})$, so
$\eaug{8a}=\eaug{8b}$ and consequently $\eaug{8a}=\eaug{8b}=0$. Further
$\chi_{5}(u)=2(\eaug{6a}-\eaug{3a})=4\eaug{6a}=\eaug{6a}\chi_{5}(1)$,
so if $\eaug{6a}\neq 0$ then $u$ is mapped
under a representation of $G$ affording $\chi_{5}$
to the identity matrix or the negative of the identity
matrix, leading to the contradiction
$\chi_{5}(1)=\chi_{5}(u^{6})=\chi_{5}(z)=-4$. Thus $\eaug{3a}=\eaug{6a}=0$.
So far, we have shown that $\eaug{12a}$ and $\eaug{12b}$ are the only
possibly non-vanishing partial augmentations of $u$.
We continue with a formal application of the Luthar--Passi method.
Let $\xi$ be a $12$-th root of unity. Then
\[ \mu(\xi,u,\chi_{6}) = \tfrac{1}{12}\big(
\text{\rm Tr}_{\sQQ(\zeta_{12})/\sQQ}(\chi_{6}(u)\xi^{-1})
+6\mu(\xi^{2},u^{2},\chi_{6})+
\text{\rm Tr}_{\sQQ(\zeta_{12}^{3})/\sQQ}(\chi_{6}(u^{3})\xi^{-3})\big). \]
Since $u^{3}$ is rationally conjugate to an element of order $4$ in $G$,
we have $\chi_{6}(u^{3})=0$.
Since $\chi_{6}(u)=\beta(\eaug{12a}-\eaug{12b})=
(\zeta_{12}^{7}-\zeta_{12}^{11})(\eaug{12a}-\eaug{12b})$, we have
\begin{align*}
&\text{\rm Tr}_{\sQQ(\zeta_{12})/\sQQ}(\chi_{6}(u)\zeta_{12}^{-7})=
6(\eaug{12a}-\eaug{12b}),\\
&\text{\rm Tr}_{\sQQ(\zeta_{12})/\sQQ}(\chi_{6}(u)\zeta_{12}^{-11})=
-6(\eaug{12a}-\eaug{12b}).
\end{align*}
Next, $\chi_{6}(u^{4})=1$ since $u^{4}$ is rationally conjugate to an element
of order $3$ in $G$, and $\chi_{6}(u^{6})=\chi_{6}(z)=-4$,
from which it is easy to see that $\mu(\xi^{2},u^{2},\chi_{6})=1$
for a primitive $12$-th root of unity $\xi$. Thus
\begin{align*}
& \mu(\zeta_{12}^{7},u,\chi_{6}) =
\tfrac{1}{2}((\eaug{12a}-\eaug{12b})+1)\geq 0, \\
& \mu(\zeta_{12}^{11},u,\chi_{6}) =
\tfrac{1}{2}(-(\eaug{12a}-\eaug{12b})+1)\geq 0,
\end{align*}
from which we obtain $\abs{\eaug{12a}-\eaug{12b}}\leq 1$. Together with
$\eaug{12a}+\eaug{12b}=1$ it follows that $\eaug{12a}=0$ or $\eaug{12b}=0$.
We have shown that all but one of the partial augmentations of $u$ vanish.

\subsection*{When $\boldsymbol{u}$ has order $\boldsymbol{8}$.}
Then the partial augmentations of $u$ which are possibly nonzero are
$\eaug{4a}$, $\eaug{4b}$, $\eaug{8a}$ and $\eaug{8b}$.
Since $\pi(u)$ has order $4$, its partial augmentations with respect to
classes of elements of order $2$ vanish and consequently
$\eaug{4a}=\eaug{4b}=0$. We have $\chi_{11}(u)=\alpha(\eaug{8a}-\eaug{8b})=
(-\zeta_{8}^{}+\zeta_{8}^{3})(\eaug{8a}-\eaug{8b})$, and this time the
Luthar--Passi method gives
\begin{align*}
& \mu(\zeta_{8}^{3},u,\chi_{11}) =
\tfrac{1}{2}((\eaug{8a}-\eaug{8b})+3)\geq 0, \\
& \mu(\zeta_{8}^{},u,\chi_{11}) =
\tfrac{1}{2}(-(\eaug{8a}-\eaug{8b})+3)\geq 0,
\end{align*}
from which we obtain $\abs{\eaug{8a}-\eaug{8b}}\leq 3$. Together with
$\eaug{8a}+\eaug{8b}=1$ it follows that
$(\eaug{8a},\eaug{8b})\in\{(1,0),(0,1),(-1,2),(2,-1)\}$.
At this point we are stuck when limiting attention to
complex characters only.

However, we may resort to $p$-modular characters.
It is natural to choose $p=5$ since $\widetilde{S}_{5}$ is a subgroup
of $\text{\rm SL}(2,25)$. This can be seen as follows.
The group $\text{\rm PSL}(2,25)$ contains $\text{\rm PGL}(2,5)$ as a
subgroup (see \cite[Kapitel~II, 8.27~Hauptsatz]{Hupp:67}) which
is isomorphic to $S_{5}$ and its pre-image in $\text{\rm SL}(2,25)$ is
isomorphic to $\widetilde{S}_{5}$ (for example, since the Sylow $2$-subgroups
of $\text{\rm SL}(2,25)$ are generalized quaternion groups).
Let $\varphi$ be the Brauer character afforded by a faithful representation
$D:\widetilde{S}_{5}\rightarrow\text{\rm SL}(2,25)$. The Brauer lift can be
chosen such that $\varphi(x)=\alpha=-\zeta_{8}^{}+\zeta_{8}^{3}$ for
an element $x$ in the conjugacy class {\tt 8a} of $G$ (since $D(x)$ has
determinant $1$). Then $x^{5}$ lies in the class {\tt 8b}, and we obtain
$\varphi(u)=\eaug{8a}\varphi(x)+\eaug{8b}\varphi(x^{5})=
(-\zeta_{8}^{}+\zeta_{8}^{3})(\eaug{8a}-\eaug{8b})$.
Since $\varphi(u)$ is the sum of two $8$-th roots of unity, it follows
that $\abs{\eaug{8a}-\eaug{8b}}\leq 1$ and consequently
$\eaug{8a}=0$ or $\eaug{8b}=0$. The proof is complete.

We remark that the choice of $p=5$ is also strongly suggested by
the general theory of cyclic blocks.
The spin characters of $\widetilde{S}_{5}$ form a single $5$-block
of defect $1$, with Brauer graph (cf.\ \cite[Theorem~4]{MoYa:88})
\begin{center}
\setlength{\unitlength}{0.8cm}
\begin{picture}(6.6,1.1)
\put(0.1,0.6){\circle{0.2}}
\put(0.2,0.6){\line(1,0){1.4}}
\put(1.7,0.6){\circle{0.2}}
\put(1.8,0.6){\line(1,0){1.4}}
\put(3.3,0.6){\circle{0.2}}
\put(3.4,0.6){\line(1,0){1.4}}
\put(4.9,0.6){\circle{0.2}}
\put(5.0,0.6){\line(1,0){1.4}}
\put(6.5,0.6){\circle{0.2}}
\put(0,0){\makebox(0.2,0.5){$\chi_{6}$}}
\put(1.6,0){\makebox(0.2,0.5){$\chi_{11}$}}
\put(3.2,0){\makebox(0.2,0.5){$\chi_{5}$}}
\put(4.8,0){\makebox(0.2,0.5){$\chi_{12}$}}
\put(6.4,0){\makebox(0.2,0.5){$\chi_{7}$}}
\put(0.8,0.7){\makebox(0.2,0.4){$\varphi_{4a}$}}
\put(2.4,0.7){\makebox(0.2,0.4){$\varphi_{2a}$}}
\put(4.0,0.7){\makebox(0.2,0.4){$\varphi_{2b}$}}
\put(5.6,0.7){\makebox(0.2,0.4){$\varphi_{4b}$}}
\end{picture}
\end{center}
Here, $\varphi=\varphi_{2a}$ and $\varphi_{2b}$ is conjugate to
$\varphi$ under the Frobenius homomorphism.

\section{The general linear group $\text{GL}(2,5)$}\label{Sec:GL25}

We set $G=\text{GL}(2,5)$. Let $z$ be a generator of $\zen{G}$, which
is a cyclic group of order $4$. The quotient
$G/\lara{x}$ is isomorphic to $S_{5}$ for which (ZC1) is known to hold.
Let $\pi$ denote the natural map $\ZZ G\rightarrow\ZZ G/\lara{z}$.

Let $u$ be a nontrivial torsion unit in $\NU{\ZZ G}$.
We will show that all but one of its partial augmentations vanish.
For that, we use part of the character table of $G$, shown in Table~\ref{T2}
in the form obtained by requiring
{\tt CharacterTable("GL25")} in GAP \cite{GAP4}, together with the natural
$2$-dimensional representation of $G$ in characteristic $5$.
In Table~\ref{T2}, the row inscribed ``in $S_{5}$'' indicates to which classes
in the quotient $S_{5}$ the listed classes of $G$ are mapped.
The classes omitted are the classes ${\tt 2a}$, ${\tt 4a}$, ${\tt 4b}$
of central $2$-elements, and the classes ${\tt 5a}$, ${\tt 20a}$,
${\tt 10a}$, ${\tt 20b}$ of elements of order divisible by $5$.

\begin{table}[h]
\[
\begin{array}{c}
\begin{array}{r|r@{\hspace*{3pt}}r@{\hspace*{3pt}}
r@{\hspace*{3pt}}r@{\hspace*{3pt}}r@{\hspace*{3pt}}
r@{\hspace*{3pt}}r@{\hspace*{3pt}}r@{\hspace*{3pt}}
r@{\hspace*{3pt}}r@{\hspace*{3pt}}r@{\hspace*{3pt}}
r@{\hspace*{3pt}}r@{\hspace*{3pt}}r@{\hspace*{3pt}}
r@{\hspace*{3pt}}r@{\hspace*{3pt}}r} \hline
\text{class} & {\tt 1a} & {\tt 4c} & {\tt 2b} &  {\tt 4d} &  {\tt 4e}
& {\tt 4f} & {\tt 4g} & {\tt 24a} & {\tt 12a} & {\tt 8a}
& {\tt 6a} & {\tt 24b} & {\tt 3a} & {\tt 8b}
& {\tt 24c} & {\tt 12b} & {\tt 24d} \rule[-0pt]{0pt}{13pt} \\
\text{in } S_{5} & {\tt\scriptstyle 1a} &
{\tt\scriptstyle 4a} & {\tt\scriptstyle 2a} & {\tt\scriptstyle 4a} &
{\tt\scriptstyle 4a} & {\tt\scriptstyle 2a} & {\tt\scriptstyle 4a} &
{\tt\scriptstyle 6a} & {\tt\scriptstyle 3a} & {\tt\scriptstyle 2b} &
{\tt\scriptstyle 3a} & {\tt\scriptstyle 6a} & {\tt\scriptstyle 3a} &
{\tt\scriptstyle 2b} & {\tt\scriptstyle 6a} & {\tt\scriptstyle 3a} &
{\tt\scriptstyle 6a}
\rule[0pt]{0pt}{-7pt} \\ \hline
\chi_{2} & 1  & i &-1 & -i & -i  &1 &  i &  i
& -1& -i & 1 & -i & 1 & i &  i & -1&  -i
\rule[0pt]{0pt}{13pt} \\
\chi_{6} & 5  & i &-1 & -i & -i  &1 &  i & -i
&  1&  i &-1 &  i &-1 &-i & -i &  1&   i\\
\chi_{16} & 4  & . & . &  . &  .  &. &  . & -i
& -1&  -2i & 1 &  i & 1 &2i & -i & -1&   i\\
\chi_{9} & 6  & \alpha & . & \bar{\alpha} &-\bar{\alpha}  &. & -\alpha
&  . &  .&  . & . &  . & . & . &  . &  .&   .\\
\chi_{14} & 6 & -\alpha & . &-\bar{\alpha} & \bar{\alpha}  &. &  \alpha
&  . &  .&  . & . &  . & . & . &  . &  .&   .\\
\chi_{15} & 4 &  . & . &  . &  .  &. &  . &  \beta & -i&  .
&-1 &-\bar{\beta} & 1 & . & -\beta &  i&  \bar{\beta}\\
\chi_{21} & 4 &  . & . &  . &  .  &. &  . &  . & 2i&  .
& 2 &  . &-2 & . &  . &  -2i&   .\\
\chi_{22} & 4 &  . & . &  . &  .  &. &  . & -\beta & -i&  .
&-1 & \bar{\beta} & 1 & . &  \beta &  i& -\bar{\beta}
\rule[-7pt]{0pt}{5pt} \\ \hline
\end{array} \\
\begin{array}{rl}
\text{Irrational entries:} &
\alpha=1+i, \\
& \beta=-\zeta+\zeta^{17} \text{ where } \zeta=\exp(2\pi i/24).
\rule[6pt]{0pt}{5pt} \\
\end{array}
\rule[0pt]{0pt}{25pt}
\end{array} \]
\rule[0pt]{0pt}{5pt}
\caption{Part of the character table of $\text{GL}(2,5)$}\label{T2}
\end{table}

The characters $\chi_{6}$ and $\chi_{16}$ have kernel $\lara{z^{2}}$.
The faithful characters $\chi_{9}$, $\chi_{14}$, $\chi_{15}$,
$\chi_{21}$ and $\chi_{22}$ of $G$ form a $5$-block of $G$,
with Brauer graph (cf.\ the theory of blocks of cyclic defect)
\begin{center}
\setlength{\unitlength}{0.8cm}
\begin{picture}(6.6,1.1)
\put(0.1,0.6){\circle{0.2}}
\put(0.2,0.6){\line(1,0){1.4}}
\put(1.7,0.6){\circle{0.2}}
\put(1.8,0.6){\line(1,0){1.4}}
\put(3.3,0.6){\circle{0.2}}
\put(3.4,0.6){\line(1,0){1.4}}
\put(4.9,0.6){\circle{0.2}}
\put(5.0,0.6){\line(1,0){1.4}}
\put(6.5,0.6){\circle{0.2}}
\put(0,0){\makebox(0.2,0.5){$\chi_{15}$}}
\put(1.6,0){\makebox(0.2,0.5){$\chi_{9}$}}
\put(3.2,0){\makebox(0.2,0.5){$\chi_{21}$}}
\put(4.8,0){\makebox(0.2,0.5){$\chi_{14}$}}
\put(6.4,0){\makebox(0.2,0.5){$\chi_{22}$}}
\put(0.8,0.7){\makebox(0.2,0.4){$\varphi_{4a}$}}
\put(2.4,0.7){\makebox(0.2,0.4){$\varphi_{2a}$}}
\put(4.0,0.7){\makebox(0.2,0.4){$\varphi_{2b}$}}
\put(5.6,0.7){\makebox(0.2,0.4){$\varphi_{4b}$}}
\end{picture}
\end{center}
Set $\varphi=\varphi_{4a}=(\chi_{15}-\chi_{9})|_{G_{5^{\prime}}}$
(restriction to $5$-regular elements). Then $\varphi$
is an irreducible $5$-modular Brauer character of $G$ of degree $2$
afforded by a natural representation $G\rightarrow\text{GL}(2,5)$.

We remark that the remaining irreducible faithful characters of $G$
form a $5$-block of $G$ which is algebraically conjugate to the one
we consider.

We write $\varepsilon_{{\tt 4c}}, \varepsilon_{{\tt 2b}},\dotsc,
\varepsilon_{{\tt 24d}}$ for the partial augmentations of $u$
at the classes listed in Table~\ref{T2}.
We assume that $u$ is not a central unit, so that its
partial augmentations at central group elements are zero.
It follows from Remark~\ref{rem1} that the order of $u$ agrees
with the order of some group element of $G$.

\subsection*{When the order of $\boldsymbol{u}$ is divisible by
$\boldsymbol{5}$.} Then $\pi(u)$ has order $5$, by Remark~\ref{rem1},
and $u$ is the product of a unit of order $5$ and a central
group element of $G$. Since there is only one class of elements of
order $5$ in $G$, the $5$-part of $u$ is rationally
conjugate to an element of $G$ (Remark~\ref{rem2}), and thus
the same is valid for $u$.

\subsection*{When $\boldsymbol{u}$ has order $\boldsymbol{2}$.}
The group $G$ has only one class of non-central elements of order $2$,
so Remark~\ref{rem2} applies.

\subsection*{When $\boldsymbol{u}$ has order $\boldsymbol{4}$.}
Then $\paug{g}{u}=0$ for a group element $g$ which is not a non-central
element of order $2$ or $4$ (Remark~\ref{rem2}). Evaluating the
Brauer character $\varphi$ at $u$ gives
\begin{equation}\label{eq2}
\varphi(u)=(\eaug{4c}-\eaug{4g})(1+i)+(\eaug{4d}-\eaug{4e})(1-i).
\end{equation}

First, suppose that $\pi(u)$ has order $2$.
Then $u^{2}=z^{2}$ (Remark~\ref{rem1}), so
$\varphi(u^{2})=-2$ and $\varphi(u)$ is the sum of two primitive
fourth roots of unity. These roots of unity are distinct since $u$
is non-central in $\ZZ G$. Thus $\varphi(u)=i+(-i)=0$ and (\ref{eq2})
gives $\eaug{4c}=\eaug{4g}$ and $\eaug{4d}=\eaug{4e}$.
Since (ZC1) holds for $S_{5}$ we have
\[ \eaug{4c}+\eaug{4g}+\eaug{4d}+\eaug{4e}=0, \qquad
\eaug{2b}+\eaug{4f}=1. \]
From that we further obtain $\eaug{4d}=-\eaug{4c}$ and
$\chi_{2}(u)=1-2\eaug{2b}+4\eaug{4c}i$. Since $|\chi_{2}(u)|=1$ it
follows that $\eaug{4c}=0$ and $\eaug{2b}\in\{0,1\}$.
Thus all but one of the partial augmentations of $u$ vanish.

Secondly, suppose that $\pi(u)$ has order $4$.
Then $\varphi(u^{2})\neq -2$. Since $\varphi(u)$ is the sum of two
distinct fourth roots of unity we have $|\varphi(u)|<2$. Thus
$\varphi(u)\in\{\pm(1+i),\pm(1-i)\}$ by (\ref{eq2}).
Since (ZC1) holds for $S_{5}$ we have
\[ \eaug{4c}+\eaug{4g}+\eaug{4d}+\eaug{4e}=1, \qquad
\eaug{2b}+\eaug{4f}=0. \]
From that and (\ref{eq2}) we further obtain that for some $a\in\ZZ$ and
$\delta_{i}\in\{0,1\}$, with exactly one $\delta_{i}$ nonzero,
$\eaug{4c}=a+\delta_{1}$, $\eaug{4g}=a+\delta_{2}$,
$\eaug{4c}=a-\delta_{3}$ and $\eaug{4g}=a-\delta_{4}$. Thus
$\chi_{2}(u)=(\delta_{1}+\delta_{2}-\delta_{3}-\delta_{4})i-2\eaug{2b}+4ai$
from which $\eaug{2b}=0$ and $a=0$ follows.
Thus all but one of the partial augmentations of $u$ vanish.

\subsection*{When $\boldsymbol{u}$ has order $\boldsymbol{8}$.}
Then $\eaug{8a}\neq 0$ or $\eaug{8b}\neq 0$ by
\cite[Corollary~4.1]{CoLi:65} (an observation sometimes attributed to
Zassenhaus). Since $S_{5}$ has no elements of order
$8$ we have $u^{4}=z^{2}$ (by Remark~\ref{rem1}).

Suppose that $\eaug{8b}=-\eaug{8a}$. Then $\chi_{16}(u)=-4\eaug{8b}i$
(remember Remark~\ref{rem2}). Since $\chi_{16}$ has degree $4$ it
follows that under a representation of $G$ affording $\chi_{16}$
the unit $u$ is mapped to a scalar multiple of the identity matrix.
Thus the image of $u$ in $\ZZ G/\lara{z^{2}}$ is a central unit of
order $4$, the kernel of $\chi_{16}$ being $\lara{z^{2}}$.
But the center of $G/\lara{z^{2}}$ has order $2$, so we have
reached a contradiction.

Hence $\eaug{8a}+\eaug{8b}\neq 0$, and since $\eaug{8a}$ and
$\eaug{8b}$ are the classes of $G$ which map onto class
$\eaug{2b}$ in $S_{5}$, in fact $\eaug{8a}+\eaug{8b}=1$.
Now $\chi_{16}(u)=2(1-2\eaug{8a})i$, and $|\chi_{16}(u)|\leq 4$
implies that $\eaug{8a}\in\{0,1\}$, so one of  $\eaug{8a}$ and
$\eaug{8b}$ vanish.

Next, we note that $\chi_{9}(u)=0$: From
$\chi_{9}(u^{4})=\chi_{9}(z^{2})=-\chi_{9}(1)$ we conclude that
$\chi_{9}(u)\in\zeta_{8}\ZZ[i]$ for a primitive $8$-th root of
unity $\zeta_{8}$ and a look at the character table shows that
$\chi_{9}(u)\in\ZZ[i]$, but definitely $\zeta_{8}\not\in\ZZ[i]$.

Since
\begin{equation}\label{eq3}
\chi_{9}(u)=(\eaug{4c}-\eaug{4g})(1+i)+(\eaug{4d}-\eaug{4e})(1-i)
\end{equation}
and $\eaug{4c}+\eaug{4g}+\eaug{4d}+\eaug{4e}=0$ it follows that
$\eaug{4c}=\eaug{4g}=-\eaug{4d}=-\eaug{4e}$. Also
$\eaug{2b}+\eaug{4f}=0$. So
$\chi_{2}(u)=(\pm 1+4\eaug{4c})i-2\eaug{2b}$ which implies
$\eaug{2b}=0$ and $\eaug{4c}=0$, and we are done.

\subsection*{When $\boldsymbol{u}$ has order $\boldsymbol{3}$.}
The group $G$ has only one class of elements of order $2$,
so Remark~\ref{rem2} applies.

\subsection*{When $\boldsymbol{u}$ has order $\boldsymbol{6}$.}
The only partial augmentations of $u$ which are possibly nonzero are
$\eaug{2b}$, $\eaug{3a}$ and $\eaug{6a}$. Since the class $\eaug{6a}$
maps in $S_{5}$ to the class of elements of order $3$
it follows that $\pi(u)$ is rationally conjugate to a group element
of order $3$ in $S_{5}$. Hence $u$ is the product of $z^{2}$ and a
unit of order $3$ (Remark~\ref{rem1}), and $u$ is rationally
conjugate to a group element.

\subsection*{When $\boldsymbol{u}$ has order $\boldsymbol{12}$.}
Only partial augmentations of $u$ taken at classes of elements of
order $2$, $4$, $3$, $6$ and $12$ are possibly nonzero.
The classes of elements of order $3$, $6$ and $12$ map in $S_{5}$
to the class of elements of order $3$. Thus $\pi(u)$ is of order $3$
and  $u$ is the product of $z$ and a unit of order $3$, so
$u$ is rationally conjugate to a group element.

\subsection*{When $\boldsymbol{u}$ has order $\boldsymbol{24}$.}
Then $\pi(u)$ is  rationally conjugate to an element of order $6$
in $S_{5}$, and so
\begin{equation}\label{eq4}
\begin{split}
& \eaug{24a}+\eaug{24b}+\eaug{24c}+\eaug{24d}=1, \\
& \eaug{12a}+\eaug{6a}+\eaug{3a}+\eaug{12b}=0, \\
& \eaug{8a}+\eaug{8b}=0, \\
& \eaug{4c}+\eaug{4d}+\eaug{4e}+\eaug{4g}=0, \\
& \eaug{2b}+\eaug{4f}=0.
\end{split}
\end{equation}

From $\chi_{9}(u^{12})=-\chi_{9}(1)$ we conclude that
$\chi_{9}(u)\in\zeta_{8}\ZZ[i,\zeta_{3}]$ for a primitive $8$-th
root of unity $\zeta_{8}$ and a primitive third root of unity
$\zeta_{3}$. A look at the character table shows
$\chi_{9}(u)\in\ZZ[i]$, so $\chi_{9}(u)=0$ as
$\zeta_{8}\not\in\ZZ[i,\zeta_{3}]$.
In the same way we argue that $\chi_{21}(u)=0$.

Thus evaluation (\ref{eq3}) of $\chi_{9}(u)$ is zero, and with
(\ref{eq4}) it follows that
$\eaug{4c}=\eaug{4g}=-\eaug{4d}=-\eaug{4e}$. Now
$(\chi_{2}+\chi_{6})(u)=-4\eaug{2a}+8\eaug{4c}i$. Since
$\chi_{2}+\chi_{6}$ has degree $6$ we conclude that $\eaug{4c}=0$.

We have
$0=\chi_{21}(u)=2(\eaug{6a}-\eaug{3a})+2i(\eaug{12a}-\eaug{12b})$,
so $\eaug{6a}=\eaug{3a}$ and $\eaug{12a}=\eaug{12b}$.
Further $\eaug{6a}=-\eaug{12a}$ from (\ref{eq4}). Thus
$\chi_{16}(u)\in -4\eaug{12a}+i\ZZ$. From
$\chi_{16}(u^{6})=-\chi_{16}(1)$ we obtain
$\chi_{16}(u)\in i\ZZ[\zeta_{3}]$. It follows that
$-4\eaug{12a}i\in\ZZ[\zeta_{3}]$ and $\eaug{12a}=0$.
Now $\chi_{2}(u)\in -2\eaug{2b}+i\ZZ$ and so $\eaug{2b}=0$.
Also $(\chi_{2}+\chi_{16})(u)=-2\eaug{2b}-6\eaug{8a}i$ and since
$\chi_{2}+\chi_{16}$ has degree $5$ we have $\eaug{8a}=0$.

Set $a=\eaug{24a}+\eaug{24c}$ and $b=\eaug{24b}+\eaug{24d}$.
Then $\chi_{2}(u)=(a-b)i$ and thus $a-b=\pm 1$. Together with
$a+b=1$ it follows that $(a,b)=(1,0)$ or $(a,b)=(0,1)$.
In the first case,
$\chi_{15}(u)=(2\eaug{24a}-1)\beta-2\eaug{24b}\bar{\beta}$,
and in the second
$\chi_{15}(u)=2\eaug{24a}\beta+(1-2\eaug{24b})\bar{\beta}$.
Using the sum formula for $\sin$ with
$\tfrac{\pi}{12}=\tfrac{\pi}{3}-\tfrac{\pi}{4}$ it is
easiest to calculate $\beta=-\sqrt{\tfrac{3}{2}}(1+i)$.
In particular, $|\beta|=\sqrt{3}$.
Since $\chi_{15}(u)$ is the sum of four roots of unity, it is
readily seen that if $\chi_{15}(u)$ assumes the first value, then
$\eaug{24b}=0$ and $\eaug{24a}\in\{0,1\}$, and
if $\chi_{15}(u)$ assumes the second value, then
$\eaug{24a}=0$ and $\eaug{24b}\in\{0,1\}$. It follows that exactly
one of $\eaug{24a}$, $\eaug{24c}$, $\eaug{24b}$ and $\eaug{24d}$
is nonzero, and we are done.
We remark that the last argument can be replaced by a
simpler ``modular'' argument:
we already know that $\chi_{15}(u)$ agrees with the value
of $\chi_{15}$ at a class of elements of order $24$
since $\chi_{9}(u)=0$ and
$\varphi=\chi_{9}-\chi_{15}$ on $5$-regular elements.



\providecommand{\bysame}{\leavevmode\hbox to3em{\hrulefill}\thinspace}
\providecommand{\MR}{\relax\ifhmode\unskip\space\fi MR }
\providecommand{\MRhref}[2]{%
  \href{http://www.ams.org/mathscinet-getitem?mr=#1}{#2}
}
\providecommand{\href}[2]{#2}

\end{document}